\documentclass[10pt]{amsart}
\usepackage{amssymb, amscd, amsmath, amsthm, epsf, epsfig, latexsym, enumerate}

\newtheorem{theorem}{Theorem}
\newtheorem{lemma}[theorem]{Lemma}
\newtheorem{cor}[theorem]{Corollary}

\begin{document}

\title{Strongly minimal $PD_4$-complexes}

\author{Jonathan A. Hillman }
\address{School of Mathematics and Statistics\\
     University of Sydney, NSW 2006\\
      Australia }

\email{jonh@maths.usyd.edu.au}

\begin{abstract}
We consider the homotopy types of $PD_4$-complexes $X$ 
with fundamental group $\pi$ such that $c.d.\pi=2$ and $\pi$ has one end. 
Let $\beta=\beta_2(\pi;\mathbb{F}_2)$ and $w=w_1(X)$.
Our main result is that (modulo two technical conditions on $(\pi,w)$)
there are at most $2^\beta$ orbits of $k$-invariants determining 
``strongly minimal" complexes (i.e., 
those with homotopy intersection pairing $\lambda_X$ trivial).
The homotopy type of a $PD_4$-complex $X$ 
with $\pi$ a $PD_2$-group is determined by $\pi$, 
$w$, $\lambda_X$ and the $v_2$-type of $X$. 
Our result also implies that Fox's 2-knot with metabelian group
is determined up to homeomorphism by its group.
\end{abstract}



\maketitle

It remains an open problem to give a homotopy classification of 
closed 4-manifolds, or more generally $PD_4$-complexes,  
in terms of standard invariants such as the fundamental
group, characteristic classes and homotopy intersection pairings.
The class of groups of cohomological dimension at most 2 seems to be
both tractable and of direct interest to geometric topology,
as it includes all surface groups, 
knot groups and the groups of many other bounded 3-manifolds. 
In our earlier papers we have shown that this case can largely be reduced
to the study of ``strongly minimal" $PD_4$-complexes
$Z$ with trivial intersection pairing on $\pi_2(Z)$.
If $X$ is a $PD_4$-complex with fundamental group $\pi$, $k_1(X)=0$
and there is a 2-connected degree-1 map $p:X\to{Z}$,
where $Z$ is strongly minimal then the homotopy type of $X$ is 
determined by $Z$ and the intersection pairing $\lambda_X$ 
on the ``surgery kernel" $K_2(p)=\mathrm{Ker}(\pi_2(p))$,
which is a finitely generated projective left $\mathbb{Z}[\pi]$-module
\cite{Hi06}.
Here we shall attempt to determine the homotopy types of 
such strongly minimal $PD_4$-complexes, 
under further hypotheses on $\pi$
and the orientation character.

The first two sections review material about generalized 
Eilenberg-Mac Lane spaces and cohomology with twisted coefficients,
the Whitehead quadratic functor and $PD_4$-complexes.
We assume thereafter that $X$ is a $PD_4$-complex,
$\pi=\pi_1(X)$ and $c.d.\pi=2$.
Such complexes have strongly minimal models $p:X\to{Z}$.
In \S3 we show that the homotopy type of $X$ is determined
by its first three homotopy groups and the second $k$-invariant
$k_2(X)\in H^4(L_\pi(\pi_2(X),2);\pi_3(X))$.

The key special cases in which the possible strongly
minimal models are when: 
\begin{enumerate}
\item $\pi\cong {F(r)}$ is a finitely generated free group;

\item $\pi=F(r)\rtimes{Z}$; or

\item $\pi$ is a $PD_2$-group.
\end{enumerate}
We review the first two cases briefly in \S4, 
and in \S5 we outline an argument 
for the case of $PD_2$-groups,
which involves cup product in integral cohomology.
(This is a model for our later work in Theorem 13.)
In Theorem 8 we show that
the homotopy type of a $PD_4$-complex $X$ 
with $\pi$ a $PD_2$-group is determined by $\pi$, 
$w=w_1(X)$, $\lambda_X$ and the $v_2$-type of $X$.
(The corresponding result was already known for $\pi$ free
and in the Spin case when $\pi$ is a $PD_2$-group.)
In \S6 we assume further that $\pi$ has one end,
and give a partial realization theorem for $k$-invariants (Theorem 9);
we do not know whether the 4-complexes we construct all satisfy
Poincar\'e duality.
In \S7 and \S8 we extend the cup product argument sketched in \S5 to
a situation involving local coefficient systems,
to establish our main result (Theorem 13).
Here we show that the number of homotopy types of minimal $PD_4$-complexes 
for $(\pi,w)$ is bounded by the order of $H^2(\pi;\mathbb{F}_2)$,
provided that $(\pi,w)$ satisfies two technical conditions.
(However we do not have an explicit invariant.)
One of these conditions fails for $\pi$ a $PD_2$-group and $w_1(\pi)$ 
or $w$ nontrivial, and thus our result is far from ideal.
Nevertheless it holds in other interesting cases, 
notably when $\pi=Z*_m$ (with $m$ even) and $w=1$. (See \S9.)
In the final section we show that if $\pi$ is the group of a fibred 
ribbon 2-knot $K$ the knot manifold $M(K)$ is determined up to 
TOP $s$-cobordism by $\pi$,
while Example 10 of Fox's ``Quick Trip Through Knot Theory" \cite{Fo62}
is determined up to TOP isotopy and reflection by its group. 

\section{generalities}

Let $X$ be a topological space with fundamental group $\pi$ and
universal covering space $\widetilde{X}$,
and let $f_{X,k}:X\to{P_k(X)}$ be the $k^{th}$ stage of the Postnikov tower 
for $X$.
We may construct $P_k(X)$ by adjoining cells of dimension at least $k+2$
to kill the higher homotopy groups of $X$.
The map $f_{X,k}$ is then given by the inclusion of $X$ into $P_k(X)$, 
and is a $(k+1)$-connected map.
In particular, $P_1(X)\simeq{K=K(\pi,1)}$ and $c_X=f_{X,1}$ is the 
classifying map for the fundamental group $\pi=\pi_1(X)$.

Let $[X;Y]_K$ be the set of homotopy classes over $K$ of maps $f:X\to{Y}$
such that $c_X=c_Yf$.
If $M$ is a left $\mathbb{Z}[\pi]$-module let $L_\pi(M,n)$ be the generalized 
Eilenberg-Mac Lane space over $K$ realizing the given action of $\pi$ on $M$.
Thus the classifying map for $L=L_\pi(M,n)$ is a principal $K(M,n)$-fibration 
with a section $\sigma:K\to{L}$.
We may view $L$ as the $ex$-$K$ loop space 
$\overline\Omega{{L_\pi}(M,n+1)}$, 
with section $\sigma$ and projection $c_L$.
Let $\mu:L\times_KL\to {L}$ be the (fibrewise) loop multiplication.
Then $\mu(id_L,\sigma{c_L})=\mu(\sigma{c_L},id_L)=id_L$ in $[L;L]_K$.
Let $\iota_{M,n}\in{H^n}(L;M)$ be the characteristic element.
The function $\theta:[X,L]_K\to{H^n}(X;M)$ given by $\theta(f)=f^*\iota_{M,n}$ 
is a isomorphism with respect to the addition on $[X,L]_K$ determined by $\mu$.
Thus $\theta(id_L)=\iota_{M,n}$,
$\theta(\sigma{c_X})=0$ and $\theta(\mu(f,f'))=\theta(f)+\theta(f')$.
(See Definition III.6.5 of \cite{Ba}.)

Let $\Gamma_W$ be the quadratic functor of J.H.C.Whitehead
and let $\gamma_A:A\to\Gamma_W(A)$ be the
universal quadratic function, for $A$ an abelian group.
The natural epimorphism from $A$ onto
$A/2A=\mathbb{F}_2\otimes{A}$ is quadratic, 
and so induces a canonical epimorphism from $\Gamma_W(A)$ to $A/2A$.
The kernel of this epimorphism is the image of the
symmetric square $A\odot{A}$.
If $A$ is a $\mathbb{Z}$-torsion-free left $\mathbb{Z}[\pi]$-module 
the sequence
\[0\to{A}\odot{A}\to\Gamma_W({A})\to{A}/2{A}\to0\]
is an exact sequence of left $\mathbb{Z}[\pi]$-modules, when
${A}\odot{A}$ and $\Gamma_W({A})$ have the diagonal left $\pi$-action.
Let $A\odot_\pi{A}=\mathbb{Z}\otimes_\pi(A\odot{A})$.

The natural map from $\Pi\odot\Pi$ to $\Gamma_W(\Pi)$ is given 
by the Whitehead product $[-,-]$, and
there is a natural exact sequence of $\mathbb{Z}[\pi]$-modules
\begin{equation}
\begin{CD}
\pi_4(X)@> hwz_4>>H_4(\widetilde{X};\mathbb{Z})\to\Gamma_W(\Pi)\to
\pi_3(X)@> hwz_3>>H_3(\widetilde{X};\mathbb{Z})\to0,
\end{CD}
\end{equation}
where $hwz_q$ is the Hurewicz homomorphism in dimension $q$.
(See Chapter 1 of \cite{Ba'}.)

Let $w:\pi\to\{\pm1\}$ be a homomorphism,
and let $\varepsilon_w:\mathbb{Z}[\pi]\to\mathbb{Z}^w$
be the $w$-twisted augmentation, given by $w$ on elements of $\pi$.
Let $I_w=\mathrm{Ker}(\varepsilon_w)$.
If $N$ is a right $\mathbb{Z}[\pi]$-module let $\overline{N}$
denote the conjugate left module determined by $g.n=w(g)n.g^{-1}$ 
for all $g\in\pi$ and $n\in{N}$.
If $M$ is a left $\mathbb{Z}[\pi]$-module 
let $M^\dagger=\overline{Hom_\pi(M,\mathbb{Z}[\pi])}$.
The higher extension modules are naturally right modules, and we
set $E^iM=\overline{Ext^i_{\mathbb{Z}[\pi]}(M,\mathbb{Z}[\pi])}$.
In particular, $E^0M=M^\dagger$ and
$E^i\mathbb{Z}=\overline{H^i(\pi;\mathbb{Z}[\pi])}$.

\begin{lemma}
Let $M$ be a $\mathbb{Z}[\pi]$-module with a finite resolution
of length $n$ and such that $E^iM=0$ for $i<n$. Then 
$Aut_\pi(M)\cong{Aut_\pi(E^nM)}$.
\end{lemma}

\begin{proof}
Since $c.d.\pi\leq2$ and $E^iM=0$ for $i<n$ the dual of a finite resolution 
for $M$ is a finite resolution for $E^nM$. 
Taking duals again recovers the original resolution, 
and so $E^nE^nM\cong{M}$.
If $f\in{Aut(M)}$ it extends to an endomorphism of the resolution
inducing an automorphism $E^nf$ of $E^nM$. 
Taking duals again gives $E^nE^nf=f$.
Thus $f\mapsto{E^nf}$ determines an isomorphism 
$Aut_\pi(M)\cong{Aut_\pi(E^nM)}$.
\end{proof}

In particular, if $\pi$ is a duality group of dimension $n$ over $\mathbb{Z}$
and $\mathcal{D}=H^n(G;\mathbb{Z}[G])$ is the dualizing module then 
$\overline{\mathcal{D}}=E^n\mathbb{Z}$ and 
$Aut_\pi(\overline{\mathcal{D}})=\{\pm1\}$.
Free groups are duality groups of dimension 1, 
while if $c.d.\pi=2$ then $\pi$ is a duality group of dimension 2
if and only if it has one end ($E^1\mathbb{Z}=0$) 
and $E^2\mathbb{Z}$ is $\mathbb{Z}$-torsion-free.
(See Chapter III of \cite{Bi}.)

\section{$PD_4$-complexes}

We assume henceforth that $X$ is a $PD_4$-complex,
with orientation character $w=w_1(X)$.
Then $\pi$ is finitely presentable and $X$ is homotopy equivalent to 
$X_o\cup_\phi{e^4}$, 
where $X_o$ is a complex of dimension at most 3 and $\phi\in\pi_3(X_o)$
\cite{Wa}.
In \cite{Hi04a} and \cite{Hi04b} we used such cellular decompositions 
to study the homotopy types of $PD_4$-complexes.
Here we shall follow \cite{Hi06} instead and rely more consistently 
on the dual Postnikov approach.

\begin{lemma}
If $\pi$ is infinite the homotopy type of $X$ is determined by $P_3(X)$.
\end{lemma}

\begin{proof}
If $X$ and $Y$ are two such $PD_4$-complexes and $h:P_3(X)\to{P_3(Y)}$ 
is a homotopy equivalence then $hf_{X,3}$ is homotopic to a map $g:X\to{Y}$.
Since $\pi$ is infinite 
$H_4(\widetilde{X};\mathbb{Z})=H_4(\widetilde{Y};\mathbb{Z})=0$.
Since $g$ is 4-connected any lift to a map 
$\tilde{g}:\widetilde{X}\to\widetilde{Y}$ is a homotopy equivalence, 
by Whitehead's Theorem, and so $g$ is a homotopy equivalence.
\end{proof}

Let $\Pi=\pi_2(X)$, with the natural left $\mathbb{Z}[\pi]$-module structure.
In Theorem 11 of \cite{Hi06}  we showed that two $PD_4$-complexes 
$X$ and $Y$ with the same strongly minimal model and with trivial 
first $k$-invariants
($k_1(X)=k_1(Y)=0$ in $H^3(\pi;\Pi)$) are homotopy equivalent if and only if
$\lambda_X\cong\lambda_Y$.
The appeal to \cite{Ru92} in the second paragraph of the proof is inadequate. 
Instead we may use the following lemma. 
(In its application we need only $P_2(X)\simeq{P_2(Y)}$, 
rather than $k_1(X)=k_1(Y)=0$).

\begin{lemma}
Let $P=P_2(X)$ and $Q=P_2(Z)$, and let $f,g:P\to{Q}$ be $2$-connected maps
such that $\pi_i(f)=\pi_i(g)$ for $i=1,2$.
Then there is a homotopy equivalence $h:P\to{P}$ such that $gh\sim{f}$.
\end{lemma}

\begin{proof}
This follows from Corollaries 2.6 and 2.7 of Chapter VIII of \cite{Ba}.
\end{proof}

\begin{lemma}
Let $Z$ be a $PD_4$-complex with a finite covering space $Z_\rho$.
Then $Z$ is strongly minimal if and only if $Z_\rho$ is strongly minimal.
\end{lemma}

\begin{proof}
Let $\pi=\pi_1(Z)$, $\rho=\pi_1(Z_\rho)$  and $\Pi=\pi_2(Z)$.
Then $\pi_2(Z_\rho)\cong\Pi|_\rho$, 
and so the lemma follows from the observations that since $[\pi:\rho]$ 
is finite $H^2(\pi;\mathbb{Z}[\pi])|_\rho\cong{H^2}(\rho;\mathbb{Z}[\rho])$
and 
$Hom_{\mathbb{Z}[\pi]}(\Pi,\mathbb{Z}[\pi])|_\rho\cong
{Hom}_{\mathbb{Z}[\rho]}(\Pi|_\rho,\mathbb{Z}[\rho])$,
as right $\mathbb{Z}[\rho]$-modules.
\end{proof}

In particular, if $v.c.d.\pi\leq2$ and $\rho$ is a torsion-free subgroup of
finite index then $c.d.\rho\leq2$, and so $\chi(Z_\rho)=2\chi(\rho)$,
by Theorem 13 of \cite{Hi06}.
Hence $[\pi:\rho]$ divides $2\chi(\rho)$,
thus bounding the order of torsion subgroups of $\pi$ if 
$\chi_{virt}(\pi)=\chi(\rho)/[\pi:\rho]\not=0$.

The next theorem gives a much stronger restriction,
under further hypotheses.

\begin{theorem}
Let $Z$ be a strongly minimal $PD_4$-complex and $\pi=\pi_1(Z)$.
Suppose that $\pi$ has one end, $v.c.d.\pi=2$ and 
$E^2\mathbb{Z}$ is free abelian.
If $\pi$ has nontrivial torsion then it is a semidirect product
$\kappa\rtimes(Z/2Z)$, where $\kappa$ is a $PD_2$-group.
\end{theorem}

\begin{proof}
Let $G$ be a torsion-free subgroup of finite index in $\pi$.
Then $H^2(\pi;\mathbb{Z}[\pi])|_G=H^2(G;\mathbb{Z}[G])$,
by Shapiro's Lemma, and so
$Aut_\pi(E^2\mathbb{Z}))\leq{Aut_G}(E^2\mathbb{Z}))=\{\pm1\}$, 
by Lemma 1.
Therefore the kernel $\kappa$ of the natural action of $\pi$ on
$\Pi=\pi_2(Z)\cong{E}^2\mathbb{Z}$ has index $[\pi:\kappa]\leq2$.
Suppose that $g\in\pi$ has prime order $p>1$.
Then $H_{s+3}(Z/pZ;\mathbb{Z})\cong{H_s}(Z/pZ;\Pi)$ for $s\geq4$, 
by Lemma 2.10 of \cite{Hi}.
In particular, $Z/pZ\cong{H_4(Z/pZ;\Pi)}$.
If $g$ acts trivially on $\Pi$ then $H_4(Z/pZ;\Pi)=0$.
Thus we may assume that $\kappa$ is torsion-free, $p=2$,
$g$ acts via multiplication by $-1$ and $\pi\cong\kappa\rtimes(Z/2Z)$.
Moreover $H_4(Z/pZ;\Pi)=\Pi/2\Pi\cong{Z/2Z}$,
and so the free abelian group $E^2\mathbb{Z}\cong\Pi$ 
must in fact be infinite cyclic.
Hence $\kappa$ is a $PD_2$-group \cite{Bo}.
\end{proof}

This result settles the question on page 67 of \cite{Hi}.

\begin{cor}
If $X$ is a $PD_4$-complex with $\pi_1(X)\cong{Z*_m\rtimes{Z/2Z}}$ 
and $m>1$ then ${\chi(X)>0}$.
\end{cor} 

\begin{proof}
Let $\rho=Z*_m$. Then $\chi(X)=\frac12\chi(X_\rho)$.
Hence $\chi(X)\geq0$, 
with equality if and only if $X_\rho$ is strongly minimal, 
by Theorem 13 of \cite{Hi06}.
In that case $X$ would be strongly minimal, by Lemma 4.
Since $\pi$ is solvable $E^2\mathbb{Z}$ is free abelian \cite{Mih}.
Therefore $X$ is not strongly minimal and so $\chi(X)>0$.
\end{proof}

\section{$c.d.\pi\leq2$}

We now assume that $c.d.\pi\leq2$.
In this case we may drop the qualification ``strongly",
for the following three notions of minimality are equivalent,
by Theorem 13 of \cite{Hi06}:
\begin{enumerate}
\item $X$ is strongly minimal;
\item  $X$ is minimal with respect
to the partial order determined by 2-connected degree-1 maps;
\item $\chi(X)=2\chi(\pi)\leq\chi(Y)$ for $Y$ any $PD_4$-complex with
$(\pi_1(Y),w_1(Y))\cong(\pi,w)$.
\end{enumerate}

We have $\Pi\cong{E^2\mathbb{Z}}\oplus{P}$, 
where $P$ is a finitely generated projective left $\mathbb{Z}[\pi]$-module,
and $X$ is minimal if and only if $P=0$.
The first $k$-invariant is trivial, since $H^3(\pi;\Pi)=0$,
and so $P_2(X)\simeq{L}=L_\pi(\Pi,2)$.
Let $\sigma$ be a section for $c_L$.
The group $E_\pi(L)$ of based homotopy classes of based self-homotopy 
equivalences of $L$ which induce the identity on $\pi$ is the group 
of units of $[L,L]_K$ with respect to composition, and 
is isomorphic to a semidirect product $H^2(\pi;\Pi)\rtimes{Aut}_\pi(\Pi)$.
(See Corollary 8.2.7 of \cite{Ba}.)

\begin{lemma}
The homotopy type of $X$ is determined by $\pi$, $\Pi$,
$\pi_3(X)$ and the orbit of $k_2(X)\in{H^4(L;\pi_3(X))}$ 
under the actions of $E_\pi(L)$ and $Aut_\pi(\pi_3(X))$.
\end{lemma}

\begin{proof}
Since these invariants determine $P_3(X)$ this follows from Lemma 2.
\end{proof}

It follows from the Whitehead sequence
(1) that $H_3(\widetilde{L};\mathbb{Z})=0$ and 
$H_4(\widetilde{L};\mathbb{Z})\cong\Gamma_W(\Pi),$
since $\widetilde{L}\simeq{K(\Pi,2)}$.
Hence the spectral sequence for the universal covering 
$p_L:\widetilde{L}\to{L}$ gives exact sequences
\[0\to{Ext^2_{\mathbb{Z}[\pi]}(\mathbb{Z},\Pi)}={H^2(\pi;\Pi)}\to
{H^2(L;\Pi)}\to
{Hom_{\mathbb{Z}[\pi]}(\Pi,\Pi)}=End_\pi(\Pi)\to0,\]
which is split by $H^2(\sigma;\Pi)$, and 
\begin{equation}
\begin{CD}
0\to{Ext^2_{\mathbb{Z}[\pi]}(\Pi,\pi_3(X))}\to{H^4(L;\pi_3(X))}@>p_L^*>>
Hom_{\mathbb{Z}[\pi]}(\Gamma_W(\Pi),\pi_3(X))\to0,
\end{CD}
\end{equation}
since $c.d.\pi\leq2$.
The right hand homomorphisms are the homomorphisms induced by $p_L$, 
in each case.
(There are similar exact sequences with coefficients
any left $\mathbb{Z}[\pi]$-module $\mathcal{A}$.)
The image of $k_2(X)$ in ${Hom}(\Gamma_W(\Pi),\pi_3(X))$ 
is a representative for $k_2(\widetilde{X})$,
and determines the middle homomorphism in the Whitehead sequence (1).
If $p_L^*k_2(X)$ is an isomorphism its orbit under the action of
$Aut_\pi(\pi_3(X))$ is unique.
If $\pi$ has one end the spectral sequence for $p_X:\widetilde{X}\to{X}$ 
gives isomorphisms 
$Ext^2_{\mathbb{Z}[\pi]}(\Pi,\mathcal{A}))\cong{H^4(X;\mathcal{A})}$
for any left $\mathbb{Z}[\pi]$-module $\mathcal{A}$,
and so $f_{X,2}$ induces splittings
$H^4(L;\mathcal{A})\cong{H^4(X;\mathcal{A})}
\oplus{H^4(\Pi,2;\mathcal{A})}^\pi$.

We wish to classify the orbits of $k$-invariants for minimal $PD_4$-complexes.
We shall first review the known cases, 
when $\pi$ is a free group or a $PD_2$-group.

\section{the known cases: free groups and semidirect products}

The cases with fundamental group a free group are well-understood.
A  minimal $PD_4$-complex $X$ with $\pi\cong {F(r)}$ free of rank $r$ 
is either $\#^r(S^1\times{S^3})$, if $w=0$,
or $\#^r(S^1\tilde\times{S^3})$, if $w\not=0$.
In \cite{Hi04a} this is established by consideration 
of the chain complex $C_*(\widetilde{X})$, 
using the good homological properties of $\mathbb{Z}[F(r)]$.
From the present point of view, if $X$ is strongly minimal $\Pi=0$,
so $L=K(\pi,1)$, $H^4(L;\pi_3(X))=0$ and $k_2(X)$ is trivial.

If $X$ is not assumed to be minimal $\Pi$ is a free $\mathbb{Z}[\pi]$-module
of rank $\chi(X)+2r-2$ and the homotopy type of $X$
is determined by the triple $(\pi,w,\lambda_X)$ \cite{Hi04a}.

The second class of groups for which the minimal models are known 
are the extensions of $Z$ by finitely generated free groups.
If $\pi=F(s)\rtimes_\alpha{Z}$ the minimal models are mapping tori 
of based self-homeomorphisms of closed 3-manifolds 
$N=\#^s(S^1\times{S^2})$ (if $w|_\nu=0$) or 
$\#^s(S^1\tilde\times{S^2})$ (if $w|_\nu\not=0$).
(See Chapter 4 of \cite{Hi}.)
Two such mapping tori are orientation-preserving homeomorphic if 
the homotopy classes of the defining self-homeomorphisms are conjugate in
the group of based self homotopy equivalences $E_0(N)$.
There is a natural representation of $Aut(F(s))$ by isotopy classes of based
homeomorphisms of $N$, 
and $E_0(N)$ is a semidirect product $D\rtimes{Aut(F(s))}$,
where $D$ is generated by Dehn twists about nonseparating 
2-spheres \cite{He77}. 
We may identify $D$ with $(Z/2Z)^s=H^1(F(s);\mathbb{F}_2)$, 
and then $E_0(N)=(Z/2Z)^s\rtimes{Aut(F(s))}$,
with the natural action of $Aut(F(s))$.

Let $f$ be a based self-homeomorphism of $N$, 
and let $M(f)$ be the mapping torus of $f$.
If $f$ has image $(d,\alpha)$ in $E_0(N)$
then $\pi=\pi_1(M(f))\cong{F(s)}\rtimes_\alpha{Z}$.
Let $\delta(f)$ be the image of $d$ in
$H^2(\pi;\mathbb{F}_2)=
H^1(F(s);\mathbb{F}_2)/(\alpha-1)H^1(F(s);\mathbb{F}_2)$.
If $g$ is another based self-homeomorphism of $N$
with image $(d',\alpha)$ and $\delta(g)=\delta(f)$
then $d-d'=(\alpha-1)(e)$ for some $e\in{D}$ and so
$(d,\alpha)$ and $(d',\alpha)$ are conjugate.
In fact this cohomology group parametrizes such homotopy types; 
see Theorem 13 for a more general result (subject to some
algebraic hypotheses).
However in this case we do not yet have explicit invariants enabling us to
decide which are the possible minimal models for a given $PD_4$-complex.
(It is a remarkable fact that if $\pi=F(s)\rtimes_\alpha{Z}$
and $\beta_1(\pi)\geq2$ then $\pi$ is such a semidirect product 
for infinitely many distinct values of $s$ \cite{Bu}.
However this does not affect our present considerations.)

It can be shown that if $N$ is a $PD_3$-complex with fundamental group $\nu$
and $\pi=\nu\rtimes_\alpha{Z}$ for some automorphism $\alpha$ 
the strongly minimal $PD_4$-complexes with fundamental group $\pi$ 
are the mapping tori of based self homotopy equivalences $h$ of $N$
which induce $\alpha$.
However if $\nu$ is not free $\alpha$ may not be nonrealizable,
and there may be $PD_4$-complexes with group $\pi$ having 
no strongly minimal model.
(See Theorem 6 of \cite{Hi06} and the subsequent construction,
for the aspherical case.)

\section{the known cases: $PD_2$-groups}

The cases with fundamental group a $PD_2$-group are 
also well understood, from a different point of view.
A minimal $PD_4$-complex $X$ with $\pi$ a $PD_2$-group is
homotopy equivalent to the total space of a $S^2$-bundle over 
a closed aspherical surface.
Thus there are two minimal models for each pair $(\pi,w)$,
distinguished by their second Wu classes.
This follows easily from the fact that the inclusion of $O(3)$ into
the monoid of self-homotopy equivalences $E(S^2)$ 
induces a bijection on components and an isomorphism on fundamental groups.
(See Lemma 5.9 of \cite{Hi}.)
However it is instructive to consider this case 
from the present point of view,
in terms of $k$-invariants, as we shall extend the following 
argument to other groups in our main result.

When $\pi$ is a $PD_2$-group and $X$ is minimal
$\Pi$ and $\Gamma_W(\Pi)$ are infinite cyclic.
The action $u:\pi\to{Aut}(\Pi)$ is given by $u(g)=w_1(\pi)(g)w(g)$
for all $g\in\pi$, by Lemma 10.3 of \cite{Hi},
while the induced action on $\Gamma_W(\Pi)$ is trivial.

Suppose first that $\pi$ acts trivially on $\Pi$. 
Then $L\simeq{K\times{CP^\infty}}$.
Fix generators $t$, $x$, $\eta$ and $z$ for 
$H^2(\pi;\mathbb{Z})$, $\Pi$, $\Gamma_W(\Pi)$ and 
$H^2(CP^\infty;\mathbb{Z})=Hom(\Pi,\mathbb{Z})$, 
respectively, such that $z(x)=1$ and $2\eta=[x,x]$.
(These groups are all infinite cyclic,
but we should be careful to distinguish the generators,
as the Whitehead product pairing of $\Pi$ with itself into 
$\Gamma_W(\Pi)$ is not the pairing given by multiplication.)
Let $t,z$ denote also the generators of ${H^2}(L;\mathbb{Z})$ induced 
by the projections to $K$ and $CP^\infty$, respectively.
Then $H^2(\pi;\Pi)$ is generated by $t\otimes{x}$,
while $H^4(L;\Gamma_W(\Pi))$ is generated by $tz\otimes\eta$ and
$z^2\otimes\eta$.
(Note that $t$ has order 2 if $w_1(\pi)\not=0$.)

The action of $[K,L]_K=[K,CP^\infty]\cong{H^2(\pi;\mathbb{Z})}$
on ${H^2}(L;\mathbb{Z})$ is generated by $t\mapsto{t}$ and $z\mapsto{z+t}$.  
The action on $H^4(L;\Gamma_W(\Pi))$ is then given by
$tz\otimes\eta\mapsto{tz}\otimes\eta$ and 
$z^2\otimes\eta\mapsto{z^2\otimes\eta+2tz\otimes\eta}$. 
There are thus two possible $E_\pi(L)$-orbits of $k$-invariants,
and each is in fact realized by the total space of
an $S^2$-bundle over the surface $K$.

If the action $u$ is nontrivial these calculations 
go through essentially unchanged with 
coefficients $\mathbb{F}_2$ instead of $\mathbb{Z}$.
There are again two possible $E_\pi(L)$-orbits of $k$-invariants,
and each is realized by an $S^2$-bundle space.
(See \S4 of \cite{Hi04b} for another account.)

In all cases the orbits of $k$-invariants correspond to
the elements of $H^2(\pi;\mathbb{F}_2)=Z/2Z$.
In fact the $k$-invariant may be detected by the Wu class.
Let $[c]_2$ denote the image of a cohomology class under reduction 
{\it mod\/} (2).
Since $k_2(X)=\pm(z^2\otimes\eta+mtz\otimes\eta)$ has image 0 in $H^4(X;\Pi)$ 
it follows that $[z]_2^2\equiv{m[tz]_2}$ in $H^4(X;\mathbb{F}_2)$.
This holds also if $\pi$ is nonorientable or the action $u$ is nontrivial,
and so $v_2(X)=m[z]_2$ and the orbit of $k_2(X)$ determine each other.

If $X$ is not assumed to be minimal its minimal models
may be determined from Theorem 7 of \cite{Hi04b}.
The enunciation of this theorem in \cite{Hi04b} is not correct;
an (implicit) quantifier over certain elements of $H^2(X;\mathbb{Z}^u)$ 
is misplaced and should be ``there is" rather than ``for all".
More precisely, where it has

\smallskip
\noindent``{\sl
and let $x\in H^2(X;\mathbb{Z}^u)$ be such that 
$(x\cup c_X^*\omega_F)[X]=1$.
Then there is a $2$-connected degree-$1$ map $h:X\to E$ such that $c_E=c_Xh$
if and only if $(c_X^*)^{-1}w_1(X)$ $=(c_E^*)^{-1}w_1(E)$, $[x]_2^2=0$ 
if $v_2(E)=0$ and $[x]_2^2=[x]_2\cup c_X^*[\omega_F]_2$ otherwise}"

\smallskip
\noindent{it should read} 

\smallskip
\noindent``{\sl
Then there is a $2$-connected degree-$1$ map $h:X\to E$ such that $c_E=c_Xh$
if and only if $(c_X^*)^{-1}w_1(X)=(c_E^*)^{-1}w_1(E)$
and there is an $x\in H^2(X;\mathbb{Z}^u)$ such that 
$(x\cup c_X^*\omega_F)[X]=1$, with $[x]_2^2=0$
if $v_2(E)=0$ and $[x]_2^2=[x]_2\cup c_X^*[\omega_F]_2$ otherwise}".

\smallskip
\noindent{The} argument is otherwise correct.
Thus if $v_2(\widetilde{X})=0$ the minimal model $Z$ is uniquely 
determined by $X$; otherwise this is not so. 
Nevertheless we have the following result. 
It shall be useful to distinguish three ``$v_2$-types" of $PD_4$-complexes:
\begin{enumerate}
\item $v_2(\widetilde{X})\not=0$ (i.e., 
$v_2(X)$ is not in the image of $H^2(\pi;\mathbb{F}_2)$ under $c_X^*$);
\item $v_2(X)=0$;
\item $v_2(X)\not=0$ but $v_2(\widetilde{X})=0$ (i.e.,
$v_2(X)$ is in $c_X^*(H^2(\pi;\mathbb{F}_2))-\{0\}$);
\end{enumerate}
(This trichotomy is due to Kreck, 
who formulated it in terms of Stiefel-Whitney
classes of the stable normal bundle of a closed 4-manifold.)

\begin{theorem}
If $\pi$ is a $PD_2$-group the homotopy type of $X$
is determined by the triple $(\pi,w,\lambda_X)$ together with its $v_2$-type.
\end{theorem}

\begin{proof}
Let $t_2$ generate $H^2(\pi;\mathbb{F}_2)$.
Then $\tau=c_X^*t_2\not=0$.
If $v_2(X)=m\tau$ and ${p:X\to{Z}}$ is
a $2$-connected degree-$1$ map then $v_2(Z)=mc_Z^*t_2$,
and so there is an unique minimal model for $X$.
Otherwise $v_2(X)\not=\tau$, and so there are 
elements $y,z\in{H^2}(X;\mathbb{F}_2)$ such that $y\cup\tau\not=y^2$ and 
$z\cup\tau\not=0$.
If $y\cup\tau=0$ and $z^2\not=0$ then $(y+z)\cup\tau\not=0$ and $(y+z)^2=0$.
Taking $x=y,z$ or $y+z$ appropriately, we have $x\cup\tau\not=0$ and $x^2=0$,
so there is a minimal model $Z$ with $v_2(Z)=0$.
In all cases the theorem now follows from the main result of \cite{Hi06}.
\end{proof}

In particular, if $C$ is a smooth projective complex curve of genus $\geq1$ 
and $X=(C\times{S^2})\#\overline{CP^2}$ is a blowup of the ruled surface 
$C\times{CP^1}=C\times{S^2}$ each of the two orientable $S^2$-bundles 
over $C$ is a minimal model for $X$.
In this case they are also minimal models in the sense 
of complex surface theory. 
(See Chapter VI.\S6 of \cite{BPV}.)
Many of the other minimal complex surfaces in the Enriques-Kodaira 
classification are aspherical, and hence strongly minimal in our sense.
However 1-connected complex surfaces are never minimal in our sense, 
since $S^4$ is the unique minimal 1-connected $PD_4$-complex
and $S^4$ has no complex structure, by a classical result of Wu.
(See Proposition IV.7.3 of \cite{BPV}.)

\section{realizing $k$-invariants}

We assume now that $\pi$ has one end.
Then $c.d.\pi=2$.
If $X$ is a $PD_4$-complex with $\pi_1(X)=\pi$
then $H_3(\widetilde{X};\mathbb{Z})=H_4(\widetilde{X};\mathbb{Z})=0$.
Hence $k_2(\widetilde{X}):\Gamma_W(\Pi)\to\pi_3(X)$ is an isomorphism, 
by the Whitehead sequence (1),
while $E_\pi(L)\cong{H^2}(\pi;\Pi)\rtimes\{\pm1\}$, 
by Corollary 8.2.7 of \cite{Ba} and Lemma 3.
Thus if $X$ is minimal its homotopy type is determined by $\pi$, 
$w$ and the orbit of $k_2(X)$.
We would like to find more explicit and accessible invariants that 
characterize such orbits.
We would also like to know which $k$-invariants give rise to $PD_4$-complexes.

Let $P(k)$ denote the Postnikov 3-stage determined by 
$k\in{H^4}(L;\mathcal{A})$.

\begin{theorem}
Let $\pi$ be a finitely presentable group with $c.d.\pi=2$ and one end, 
and let $w:\pi\to\{\pm1\}$ be a homomorphism.
Let $\Pi=E^2\mathbb{Z}$ and let $k\in{H^4}(L;\Gamma_W(\Pi))$.
Then 
\begin{enumerate}
\item
There is a finitely dominated $4$-complex $Y$ with 
$H_3(\widetilde{Y};\mathbb{Z})=H_4(\widetilde{Y};\mathbb{Z})=0$ 
and Postnikov $3$-stage $P(k)$ if and only if 
$p_L^*k$ is an isomorphism and $P(k)$ has finite $3$-skeleton.
These conditions determine the homotopy type of $Y$.

\item
If $\pi$ is of type $FF$ we may assume that $Y$ is a finite complex.

\item $H_4(Y;\mathbb{Z}^w)\cong\mathbb{Z}$ and 
there are isomorphisms 
$\overline{H^p(Y;\mathbb{Z}[\pi])}\cong{H_{4-p}}(Y;\mathbb{Z}[\pi])$ 
induced by cap product with a generator $[Y]$, for $p\not=2$.
\end{enumerate}
\end{theorem}

\begin{proof}
Let $Y$ be a finitely dominated 4-complex with 
$H_3(\widetilde{Y};\mathbb{Z})=H_4(\widetilde{Y};\mathbb{Z})=0$ 
and Postnikov $3$-stage $P(k)$.
Since $Y$ is finitely dominated it is homotopy equivalent 
to a 4-complex with finite 3-skeleton, and since 
$P(k)\simeq{Y}\cup{e^{q\geq5}}$ may be constructed by adjoining cells of 
dimension at least 5 we may assume that $P(k)$ has finite 3-skeleton.
The homomorphism $p_L^*k$ is an isomorphism,
by the exactness of the Whitehead sequence (1).

Suppose now that $p_L^*k$ is an isomorphism and $P(k)$ has finite $3$-skeleton.
Let $P=P(k)^{[4]}$ and let $C_*=C_*(\widetilde{P})$ be the equivariant 
cellular chain complex for $\widetilde{P}$.
Then $C_q$ is finitely generated for $q\leq3$.
Let $B_q\leq{Z_q}\leq{C_q}$ be the submodules
of $q$-boundaries and $q$-cycles, respectively.
Clearly $H_1(C_*)=0$ and $H_2(C_*)\cong\Pi$,
while $H_3(C_*)=0$, since $p_L^*k$ is an isomorphism.
Hence there are exact sequences
\[0\to{B_1}\to{C_1}\to{C_0}\to\mathbb{Z}\to0\]
and
\[\quad  0\to{B_3}\to{C_3}\to{Z_2}\to\Pi\to0.\]
Schanuel's Lemma implies that $B_1$ is projective, since $c.d.\pi=2$.
Hence $C_2\cong{B_1}\oplus{Z_2}$ and so $Z_2$ is 
finitely generated and projective.
It then follows that $B_3$ is also finitely generated and projective,
and so $C_4\cong{B_3}\oplus{Z_4}$.
Thus $H_4(C_*)=Z_4$ is a projective 
direct summand of $C_4$.

After replacing $P$ by $P\vee{W}$,
where $W$ is a wedge of copies of $S^3$,
if necessary, we may assume that $Z_4=H_4(P;\mathbb{Z}[\pi])$ is free.
Since $\Gamma_W(\Pi)\cong\pi_3(P)$ the Hurewicz homomorphism 
from $\pi_4(P)$ to $H_4(P;\mathbb{Z}[\pi])$ is onto.
(See Chapter I\S3 of \cite{Ba'}.)
We may then attach 5-cells along maps
representing a basis to obtain a countable 5-complex $Q$
with 3-skeleton $Q^{[3]}=P(k)^{[3]}$ and
with $H_q(\widetilde{Q};\mathbb{Z})=0$ for $q\geq3$.
The inclusion of $P$ into $P(k)$ extends to a 4-connected map 
from $Q$ to $P(k)$.
Now $C_*(\widetilde{Q})$ is chain homotopy equivalent to
the complex obtained from $C_*$ by replacing $C_4$ by $B_3$,
which is a finite projective chain complex.
It follows from the finiteness conditions of Wall that $Q$ 
is homotopy equivalent to a finitely dominated complex $Y$ 
of dimension $\leq4$ \cite{Wa66}.
The homotopy type of $Y$ is uniquely determined by the data, as in Lemma 1. 

If $\pi$ is of type $FF$ then $B_1$ is stably free, by Schanuel's Lemma.
Hence $Z_2$ is also stably free.
Since dualizing a finite free resolution of $\mathbb{Z}$
gives a finite free resolution of $\Pi=E^2\mathbb{Z}$ we see in turn that
$B_3$ must be stably free, and so $C_*(\widetilde{Y})$ is chain homotopy
equivalent to a finite free complex.
Hence $Y$ is homotopy equivalent to a finite 4-complex \cite{Wa66}.

Let $D_*$ and $E_*$ be the subcomplexes of $C_*$ corresponding to the above
projective resolutions of $\mathbb{Z}$ and $\Pi$.
(Thus $D_0=C_0$, $D_1=C_1$, $D_2=B_1$ and $D_q=0$ for $q\not=0,1,2$,
while $E_2=Z_2$, $E_3=C_3$, $E_4=B_3$ and $E_r=0$ for $r\not=2,3,4$.)
Then $C_*(\widetilde{Y})\simeq{D_*}\oplus{E}_*$.
(The splitting reflects the fact that $c_Y$ is a retraction, 
since $k_1(Y)=0$.)
Clearly $H^p(Y;\mathbb{Z}[\pi])=H_{4-p}(Y;\mathbb{Z}[\pi])=0$ 
if $p\not=2$ or 4, while $H^4(Y;\mathbb{Z}[\pi])=E^2\Pi\cong\mathbb{Z}$
and $H_4(Y;\mathbb{Z}^w)=Tor_2(\mathbb{Z}^w;\Pi)\cong
\mathbb{Z}^w\otimes_\pi\mathbb{Z}[\pi]\cong\mathbb{Z}.$
The homomorphism
$\varepsilon_{w\#}:H^4(Y;\mathbb{Z}[\pi])\to{H^4}(Y;\mathbb{Z}^w)$
induced by $\varepsilon_w$ is surjective,
since $Y$ is 4-dimensional, and therefore is an isomorphism.
Hence $-\cap[Y]$ induces isomorphisms in degrees other than 2.
\end{proof}

Since $\overline{H^2(Y;\mathbb{Z}[\pi])}\cong{E^2\mathbb{Z}}$,
$H_2(Y;\mathbb{Z}[\pi])=\Pi$ and 
$Hom_\pi(E^2\mathbb{Z},\Pi)\cong{End_\pi}(E^2\mathbb{Z})$
$=\mathbb{Z}$,
cap product with $[Y]$ in degree 2 is determined by an integer,
and $Y$ is a $PD_4$-complex if and only if this integer is $\pm1$.
The obvious question is: what is this integer?
Is it always $\pm1$?
The complex $C_*$ is clearly chain homotopy equivalent to its dual,
but is the chain homotopy equivalence given by slant product with $[Y]$?

There remains also the question of characterizing the $k$-invariants
corresponding to Postnikov 3-stages with finite 3-skeleton.

If $\pi$ is either a semidirect product $F(s)\rtimes{Z}$
or the fundamental group of a Haken $3$-manifold $M$
then $\widetilde{K}_0(\mathbb{Z}[\pi])=0$,
i.e., projective $\mathbb{Z}[\pi]$-modules are stably free \cite{Wd78}.
(This is not yet known for all torsion-free one relator groups.)
In such cases finitely dominated complexes are homotopy finite.

\section{a lemma on cup products}

In our main result (Theorem 13) we shall use a ``cup-product" argument to
relate cohomology in degrees 2 and 4.
Let $G$ be a group and let $\Gamma=\mathbb{Z}[G]$.
Let $C_*$ and $D_*$ be chain complexes of left $\Gamma$-modules
and $\mathcal{A}$ and $\mathcal{B}$ left $\Gamma$-modules.
Using the diagonal homomorphism from $G$ to $G\times{G}$
we may define {\it internal products\/}
\[H^p(Hom_\Gamma(C_*,\mathcal{A}))\otimes
{H^q}(Hom_\Gamma(D_*,\mathcal{B}))\to
{H^{p+q}}(Hom_\Gamma(C_*\otimes{D_*},
\mathcal{A}\otimes\mathcal{B}))\]
where the tensor products of $\Gamma$-modules
are taken over $\mathbb{Z}$ and have the diagonal $G$-action.
(See Chapter XI.\S4 of \cite{CE}.)
If $C_*$ and $D_*$ are resolutions of $\mathcal{C}$ and $\mathcal{D}$,
respectively, we get pairings
\[ Ext^p_\Gamma(\mathcal{C},\mathcal{A})\otimes
{Ext^q_\Gamma}(\mathcal{D},\mathcal{B})\to
{Ext^{p+q}_\Gamma}
(\mathcal{C}\otimes\mathcal{D},\mathcal{A}\otimes\mathcal{B}).\]
When $\mathcal{A}=\mathcal{B}=\mathcal{D}=\Pi$, $\mathcal{C}=\mathbb{Z}$
and $q=0$ we get pairings
\[H^p(\pi;\Pi)\otimes{End}_\pi(\Pi)\to
{Ext}^p_{\mathbb{Z}[\pi]}(\Pi,\Pi\otimes\Pi).\]
If instead $C_*=D_*=C_*(\widetilde{S})$ for some space $S$ 
with $\pi_1(S)\cong{G}$ composing with an equivariant diagonal approximation gives pairings
\[ H^p(S;\mathcal{A})\otimes{H^q}(S;\mathcal{B})\to
{H^{p+q}}(S;\mathcal{A}\otimes\mathcal{B}).\]
These pairings are compatible with the universal coefficient spectral sequences
$Ext^q_\Gamma(H_p(C_*),\mathcal{A})\Rightarrow 
{H^{p+q}}(C^*;\mathcal{A})=H^{p+q}(Hom_\Gamma(C_*,\mathcal{A}))$, etc.
We shall call these pairings ``cup products",
and use the symbol $\cup$ to express their values.

We wish to show that if $c.d.\pi=2$ and $\pi$ has one end 
the homomorphism
$c^2_{\pi,w}:H^2(\pi;\Pi)\to{Ext^2_{\mathbb{Z}[\pi]}(\Pi,\Pi\otimes\Pi)}$
given by cup product with $id_{\Pi}$ is an isomorphism.
The next lemma shows that these groups are isomorphic; 
we state it in greater generality than we need, 
in order to clarify the hypotheses on the group.

\begin{lemma}
Let $G$ be a group for which the augmentation (left) module $\mathbb{Z}$ 
has a finite projective resolution $P_*$ of length $n$,
and such that $H^j(G;\Gamma)=0$ for $j<n$. 
Let $\mathcal{D}=H^n(G;\Gamma)$, $w:G\to\{\pm1\}$ be a homomorphism
and $\mathcal{B}$ be a left $\Gamma$-module.
Then there are natural isomorphisms
\begin{enumerate}
\item 
$h_{\mathcal{B}}:H^n(G;\mathcal{B})\to\mathcal{D}\otimes_G\mathcal{B}$; and

\item 
$e_{\mathcal{B}}:Ext^n_\Gamma(\overline{\mathcal{D}},\mathcal{B})\to
\mathbb{Z}^w\otimes_G\mathcal{B}=\mathcal{B}/I_w\mathcal{B}$.
\end{enumerate}
Hence $\theta_{\mathcal{B}}=
e_{\overline{\mathcal{D}}\otimes\mathcal{B}}^{-1}h_{\mathcal{B}}:
H^n(G;\mathcal{B})\cong 
Ext^n_\Gamma(\overline{\mathcal{D}},
\overline{\mathcal{D}}\otimes\mathcal{B})$ is an isomorphism;
\end{lemma}

\begin{proof}
We may assume that $P_0=\Gamma$.
Let $Q_j=Hom_\Gamma(P_{n-j},\Gamma)$
and $\partial_i^Q=Hom_\Gamma(\partial^P_{n-j},\Gamma)$.
This gives a resolution $Q_*$ for $\mathcal{D}$ 
by finitely generated projective right modules,
with $Q_n=\Gamma$.
Let $\eta:{Q}_0\to{\mathcal{D}}$ 
be the canonical epimorphism.
Tensoring $Q_*$ with $\mathcal{B}$ gives (1).
Conjugating and applying $Hom_\Gamma(-,\mathcal{B})$ gives (2).
Since we may identify ${\mathcal{D}\otimes_G\mathcal{B}}$ with
$\mathbb{Z}^w\otimes_G(\overline{\mathcal{D}}\otimes\mathcal{B})$,
composition gives an isomorphism
$\theta_{\mathcal{B}}=
e_{\overline{\mathcal{D}}\otimes\mathcal{B}}^{-1}h_{\mathcal{B}}:
H^n(G;\mathcal{B})\cong 
Ext^n_\Gamma(\overline{\mathcal{D}},
\overline{\mathcal{D}}\otimes\mathcal{B})$.
\end{proof}

If $\mathcal{D}$ is $\mathbb{Z}$-torsion free then $G$ is a
duality group of dimension $n$, with dualizing module $\mathcal{D}$.
(See \cite{Bi}.) 
It is not known whether all the groups considered in the lemma
are duality groups, even when $n=2$.

Let $A:Q_0\otimes_G\overline{\mathcal{D}}\to
{Hom}_\Gamma(P_n,\overline{\mathcal{D}})$ 
be the homomorphism given by
${A(q\otimes_G\delta)(p)}\!=q(p)\delta$ 
for all $p\in{P_n}$, $q\in{Q_0}$ and
$\delta\in\overline{\mathcal{D}}$,
and let $[\xi]\in H^n(G;\overline{\mathcal{D}})$
be the image of $\xi\in{Hom}_\Gamma(P_n,\overline{\mathcal{D}})$.
If $\xi=A(q\otimes_G\delta)$ then 
$h_{\overline{\mathcal{D}}}([\xi])=\eta(q)\otimes\delta$ and
$\xi\otimes\eta:
P_n\otimes\overline{Q}_0\to\overline{\mathcal{D}}\otimes\overline{\mathcal{D}}$
represents $[\xi]\cup{id_{\overline{\mathcal{D}}}}$
in $Ext^n_\Gamma(\overline{\mathcal{D}},
\overline{\mathcal{D}}\otimes\overline{\mathcal{D}})$.
There is a chain homotopy equivalence 
$j_*:\overline{Q}_*\to{P_*\otimes\overline{Q}_*}$,
since $P_*$ is a resolution of $\mathbb{Z}$.
Given such a chain homotopy equivalence, 
$e_{\overline{\mathcal{D}}\otimes\overline{\mathcal{D}}}
([\xi]\cup{id_{\overline{\mathcal{D}}}})$ is the image of
$(\xi\otimes\eta)(j_n(1^*))$,
where $1^*$ is the canonical generator of $\overline{Q}_n$, 
defined by $1^*(1)=1$.

Let $\tau$ be the ($\mathbb{Z}$-linear) involution of 
$H^n(G;\overline{\mathcal{D}})$ given by
$\tau(h_{\overline{\mathcal{D}}}^{-1}(\rho\otimes_G\alpha))=
h_{\overline{\mathcal{D}}}^{-1}(\alpha\otimes_G\rho))$.
If $G$ is a $PD_n$-group then $H^n(G;\overline{\mathcal{D}})\cong {Z}$
(if $w=w_1(\pi)$) or $Z/2Z$ (otherwise), and so $\tau$ is the identity.

Suppose now that  $c.d.G=2$ and $G$ has one end (i.e., $n=2$).
In order to make explicit calculations we shall assume there is a finite 
2-dimensional $K(G,1)$-complex with corresponding presentation
$\langle{X}\mid{R}\rangle$.
Then the free differential calculus gives a free resolution
\[0\to{P_2}=\Gamma\langle{p_r^2;r\in{R}}\rangle\to
{P_1}=\Gamma\langle{p_x^1;x\in{X}}\rangle\to{P_0}=\Gamma\to\mathbb{Z}\to0\]
in which 
$\partial{p_r^2}=\Sigma_{x\in{X}}r_xp_x^1$,
where $r_{x} =\frac{\partial{r}}{\partial{x}}$
and $\partial{p_x^1}=x-1$, for $r\in{R}$ and $x\in{X}$.
Let $\{q_x^1\}$ and $\{q_r^0\}$ be the dual bases for $\overline{Q}_1$
and $\overline{Q}_0$, respectively.
(Thus $q_x^1(p_y^1)=1$ if $x=y$ and 0 otherwise,
and $q_r^0(p_s^2)=1$ if $r=s$ and 0 otherwise.)
For simplicity of notation we shall write $\bar{g}=w(g)g^{-1}$ for $g\in{G}$.
Then $\partial 1^*=\Sigma_{x\in{X}}(\overline{x}-1)q_x^1$
and $\partial{q_x^1}=\Sigma_{r\in{R}}\overline{r_x}q_r^0$.
We may write $\overline{r_x}=\Sigma_k{e_{rxk}}r_{xk}$,
where $e_{rxk}=\pm1$ and $r_{xk}\in{G}$.
Then 
$r_{xk}-1=\partial(\Sigma_{y\in{X}}\frac{\partial{r_{xk}}}{\partial{y}}p_y^1)$.
Define $j_*$ in degrees 0 and 1 by setting 
\[j_0(q_r^0)=1\otimes{q_r^0}\quad\mathrm{ for }\quad{r\in{R}}\quad\mathrm
{and}\]
\[j_1(q_x^1)=1\otimes{q_x^1}+
\Sigma_{r,k,y}e_{rxk}(\frac{\partial{r_{xk}}}{\partial{y}}p_y^1\otimes
{r_{xk}q_r^0)}
\quad\mathrm{ for }\quad{x\in{X}}.\]
At this point we must specialize further.
We shall give several simple examples,
where we have managed to determine $j_2(1^*)$.
(We do not need formulae for the higher degree terms.)
The evidence suggests that if $w$ is trivial we should expect
\[j_2(1^*)=1\otimes1^*-\Sigma_{x\in{X}}x^{-1}(p_x^1\otimes{q_x^1})-\Psi\]
where $\Psi=\Sigma_{r\in{R}}u_r(p_r^2\otimes{q_r^0})$ 
with $u_r$ the inverse of a segment of $r$ and such that
\[
\partial\Psi=1\otimes\partial1^*-
\Sigma_{x\in{X}}{x}^{-1}((x-1)\otimes{q_x^1})
+\Sigma_{x\in{X}}{x}^{-1}(p_x^1\otimes\Sigma_{r\in{R}}\overline{r_x}q_r^0)
-j_1(\partial1^*)\]
\[
=\Sigma_{x,r,k}e_{rxk}[x^{-1}
((p_x^1-\Sigma_y\frac{\partial{r_{xk}}}{\partial{y}}p_y^1)\otimes{r_{xk}}q_r^0)
+(\Sigma_y\frac{\partial{r_{xk}}}{\partial{y}}p_y^1)\otimes{r_{xk}}q_r^0)].
\]

\smallskip
\noindent{\bf Examples.}
\begin{enumerate}
\item
Let $G=F(X)\times{Z}$, 
with presentation $\langle{t,X}\mid{txt^{-1}x^{-1}~\forall{x}\in{X}}\rangle$.
Then we may take
\[j_2(1^*)=1\otimes1^*-t^{-1}(p_t^1\otimes{q_t^1})-\Sigma_{x\in{X}}x^{-1}(p_x^1\otimes{q_x^1})
-\Sigma_{x\in{X}}x^{-1}t^{-1}(p_x^2\otimes{q_x^0}).\]

\medskip
\item
Let $G$ be the group with presentation $\langle{a,b}\mid{a^mb^{-n}}\rangle$.
Then we may take
\[j_2(1^*)=1\otimes1^*-a^{-1}(p_a^1\otimes{q_a^1})-b^{-1}(p_b^1\otimes{q_b^1})
-a^{-m}(p^2\otimes{q^0}).\]

\medskip
\item
Let $G$ be the orientable $PD_2$-group of genus 2, with presentation

$\langle{a,b,c,d}\mid{aba^{-1}b^{-1}cdc^{-1}d^{-1}}\rangle$.
Then we may take
\[j_2(1^*)=1\otimes1^*-\Sigma_{x\in{X}}
\overline{x}(p_x^1\otimes{q_x^1})-bab^{-1}a^{-1}(p^2\otimes{q^0}).\]

\medskip
\item
Let $G=Z*_m$ be the group with presentation 
$\langle{a,t}\mid{tat^{-1}a^{-m}}\rangle$, for $m\not=0$.
Then we may take
\[j_2(1^*)=1\otimes1^*-
a^{-1}(p_a^1\otimes{q_a^1})-t^{-1}(p_t^1\otimes{q_t^1})
-a^{-1}t^{-1}(p^2\otimes{q^0}).\]

\end{enumerate}

In each of these cases we find that 
$[\xi]\cup{id_{\overline{\mathcal{D}}}}=
-\theta_{\overline{\mathcal{D}}}(\tau([\xi]))$ 
for $\xi\in{H^2}(\pi;\overline{\mathcal{D}})$.
Similar formulae apply for $n\leq1$,
i.e., for free groups of finite rank $r\geq0$.

If $H$ is a subgroup of finite index in $G$ 
and $\mathcal{A}$ is a left $\mathbb{Z}[G]$-module 
then Shapiro's Lemma gives isomorphisms 
$H^n(G;\mathcal{A})\cong{H^n}(H;\mathcal{A}|_H)$.
Thus if $G$ satisfies the hypotheses of Lemma 10
$\overline{\mathcal{D}}|_H$ is the corresponding module for $H$.
Further applications of Shapiro's Lemma imply that
cup product with $id_{\overline{\mathcal{D}}}$ 
is an isomorphism for $(G,w)$ if and only if it is so for $(H,w|_H)$.
In particular, Examples (1)--(3) imply that $c^2_{\pi,w}$ is 
an isomorphism for all torus knot groups and $PD_2$-groups,
and all orientation characters $w$. 

\section{the action of $E_\pi(L)$}

In this section we shall attempt to study the action of $E_\pi(L)$ 
on the set of possible $k$-invariants for a minimal $PD_4$-complex 
by extending the argument sketched above for the case of $PD_2$-groups.
We believe that the restrictions we impose here on the pair $(\pi,w)$
shall ultimately be seen to be unnecessary.

Our argument shall involve relating the algebraic and homotopical 
(obstruction-theoretic) interpretations of cohomology classes.
We shall use the following special case of 
a result of Tsukiyama \cite{Tsu}; 
we give only the part that we need below.

\begin{lemma}
There is an exact sequence
$0\to{H^2}(\pi;\Pi)\to{E_\pi(L)}\to{Aut_\pi(\Pi)}\to0.$
\end{lemma}

\begin{proof}
Let $\theta:[K,L]_K\to{H^2}(\pi;\Pi)$ be the isomorphism given by
$\theta(s)=s^*\iota_{\Pi,2}$, 
and let $\theta^{-1}(\phi)=s_\phi$ for $\phi\in{H^2}(\pi;\Pi)$.
Then $s_\phi$ is a homotopy class of sections of $c_L$, 
$s_0=\sigma$ and $s_{\phi+\psi}=\mu(s_\phi,s_\psi)$,
while $\phi=s_\phi^*\iota_{\Pi,2}$.
(Recall that $\mu:L\times_KL\to {L}$ is the fibrewise loop multiplication.)

Let $h_\phi=\mu(s_\phi{c_L},id_L)$.
Then $c_Lh_\phi=c_L$ and so $h_\phi\in[L;L]_K$.
Clearly $h_0=\mu(\sigma{c_L},id_L)=id_L$
and $h_\phi^*\iota_{\Pi,2}=\iota_{\Pi,2}+c_L^*\phi\in{H^2}(L;\Pi)$.
We also see that
\[h_{\phi+\psi}=\mu(\mu(s_\phi,s_\psi){c_L},id_L)=
\mu(\mu(s_\phi{c_L},s_\psi{c_L}),id_L)=
\mu(s_\phi{c_L},\mu(s_\psi{c_L},id_L))\]
(by homotopy associativity of $\mu$) and so
\[h_{\phi+\psi}=\mu(s_\phi{c_L},h_\psi)=
\mu(s_\phi{c_L}h_\psi,h_\psi)=h_\phi{h_\psi}.\]
Therefore $h_\phi$ is a homotopy equivalence for all 
$\phi\in{H^2}(\pi;\Pi)$, 
and $\phi\mapsto{h_\phi}$ defines a homomorphism from
${H^2}(\pi;\Pi)$ to ${E_\pi(L)}$.

The lift of $h_\phi$ to the universal cover $\widetilde{L}$ is 
(non-equivariantly) homotopic to the identity, 
since the lift of $c_L$ is (non-equivariantly) homotopic to a constant map.
Therefore $h_\phi$ acts as the identity on $\Pi$.
The homomorphism $h:\phi\mapsto{h_\phi}$ is in fact an isomorphism 
onto the kernel of the action of ${E_\pi(L)}$ on $\Pi=\pi_2(L)$ \cite{Tsu}.
\end{proof}

Note also that we may view elements of $[K,L]_K$ (etc.) as
$\pi$-equivariant homotopy classes of $\pi$-equivariant maps
from $\widetilde{K}$ to $\widetilde{L}$.

\begin{lemma}
There is an exact sequence 
$\Pi\odot_\pi\Pi\to
\mathbb{Z}\otimes_\pi\Gamma_W(\Pi)\to{H^2}(\pi;\mathbb{F}_2)\to0.$
If $\Pi\odot_\pi\Pi$ is $2$-torsion-free this sequence is short exact.
\end{lemma}

\begin{proof}
Since $\pi$ is finitely presentable $\Pi$ is $\mathbb{Z}$-torsion-free 
\cite{GM}, 
and so the natural map from $\Pi\odot\Pi$ to $\Gamma_W(\Pi)$ is injective.
Applying $\mathbb{Z}\otimes_\pi-$ to the exact sequence 
\[0\to\Pi\odot\Pi\to\Gamma_W(\Pi)\to\Pi/2\Pi\to0\]
gives the above sequence,
since $\mathbb{Z}\otimes_\pi\Pi/2\Pi\cong\Pi/(2,I_w)\Pi\cong
{H^2}(\pi;\mathbb{F}_2)$.
The kernel on the left in this sequence
is the image of the 2-torsion group
$Tor_1^{\mathbb{Z}[\pi]}(\mathbb{Z},\Pi/2\Pi)$.
\end{proof}

\begin{theorem}
Let $\pi$ be a finitely presentable group such that $c.d.\pi=2$
and $\pi$ has one end.
Let $\Pi=E^2\mathbb{Z}$ and $\beta=\beta_2(\pi;\mathbb{F}_2)$.
Assume that $c^2_{\pi,w}$ is surjective and 
$\mathbb{Z}^w\otimes_\pi\Gamma_W(\Pi)$ is $2$-torsion-free.
Then there are at most $2^\beta$ orbits of $k$-invariants of minimal 
$PD_4$-complexes with Postnikov $2$-stage $L$ under the actions of $E_\pi(L)$ 
and $Aut_\pi(\Gamma_W(\Pi))$.
\end{theorem}

\begin{proof}
Let $\phi\in{H^2}(\pi;\Pi)$ and let $s_\phi\in[K,L]_K$ and $h_\phi\in[L,L]_K$ 
be as defined in Lemma 11.
Let $M={L_\pi}(\Pi,3)$ and let $\overline\Omega:{[M,M]_K}\to[L,L]_K$
be the loop map.
Since $c.d.\pi=2$ we have $[M,M]_K\cong{H^3}(M;\Pi)=End_\pi(\Pi)$.
Let $g\in[M,M]_K$ have image $[g]=\pi_3(g)\in{End_\pi(\Pi)}$
and let $f=\overline\Omega{g}$.
Then $\omega([g])=f^*\iota_{\Pi,2}$ defines a homomorphism 
$\omega:End_\pi(\Pi)\to{H^2}(L;\Pi)$ such that $p_L^*\omega([g])=[g]$ 
for all $[g]\in{End_\pi(\Pi)}$.
Moreover $f\mu=\mu(f,f)$,
since $f=\overline\Omega{g}$, and so $fh_\phi=\mu(fs_\phi{c_L},f)$.
Hence $h_\phi^*\xi=\xi+c_L^*s_\phi^*\xi$ for 
$\xi=\omega([g])=f^*\iota_{\Pi,2}$.

Naturality of the isomorphisms 
$H^2(X;\mathcal{A})\cong[X,L_\pi(\mathcal{A},2)]_K$ for $X$ a space over $K$
and $\mathcal{A}$ a left $\mathbb{Z}[\pi]$-module implies that
\[s_\phi^*\omega([g])=[g]_\#s_\phi^*\iota_{\Pi,2}=[g]_\#\phi\]
for all $\phi\in{H^2}(\pi;\Pi)$ and $g\in[M,M]_K$.
(See Chapter 5.\S4 of \cite{Ba0}.)
If $u\in{H^2}(\pi;\mathcal{A})$ then $h_\phi^*c_L^*(u)=c_L^*(u)$,
since $c_Lh_\phi=c_L$.
The homomorphism induced on the quotient 
$H^2(L;\mathcal{A})/c_L^*H^2(\pi;\mathcal{A})\cong
{Hom_{\mathbb{Z}[\pi]}(\Pi,\mathcal{A})}$
by $h_\phi$ is also the identity, 
since the lifts of $h_\phi$ are (non-equivariantly) homotopic 
to the identity in $\widetilde{L}$.

Taking $\mathcal{A}=\Pi$ we obtain a homomorphism 
$\delta_\phi:End_\pi(\Pi)\to{H^2}(\pi;\Pi)$ such that
$h_\phi^*(\xi)=\xi+c_L^*\delta_\phi(p_L^*\xi)$ for all $\xi\in{H^2}(L;\Pi)$.
Since $p_L^*\delta_\phi=0$ and $h_{\phi+\psi}=h_\phi{h_\psi}$ 
it follows that $\delta_\phi$ is additive as a function of $\phi$.
If $g\in[M,M]_K$ and 
$\phi=\rho\otimes_\pi\alpha\in{H^2}(\pi;\mathbb{Z}[\pi])\otimes_\pi\Pi$ 
then 
\[\delta_\phi([g])=\delta_\phi(p_L^*\omega([g]))=
s_\phi^*\omega[g]=\rho\otimes_\pi[g](\alpha).\]
The automorphism of $H^4(L;\mathcal{A})$ induced by $h_\phi$ 
preserves the subgroup $Ext^2_{\mathbb{Z}[\pi]}(\Pi,\mathcal{A})$ 
and induces the identity on the quotient 
$Hom_\pi(\Gamma_W(\Pi),\mathcal{A})$. 
Taking $\mathcal{A}=\Gamma_W(\Pi)$ we obtain a homomorphism 
$f_\phi=h_\phi^*-id$ from $H^4(L;\Gamma_W(\Pi))$
to $Ext^2_{\mathbb{Z}[\pi]}(\Pi,\Gamma_W(\Pi))$.

When $S=L$, $\mathcal{A}=\mathcal{B}=\Pi$, and $p=q=2$ 
the construction of \S7 gives a cup product pairing of $H^2(L;\Pi)$ 
with itself with values in ${H^4(L;\Pi\otimes\Pi)}$.
Since $c.d.\pi=2$ this pairing is trivial on the image 
of $H^2(\pi;\Pi)\otimes{H^2(\pi;\Pi)}$.
The maps $c_L$ and $\sigma$ induce a splitting
$H^2(L;\Pi)\cong{H^2(\pi;\Pi)}\oplus{End_\pi(\Pi)}$,
and this pairing restricts to the cup product pairing of 
$H^2(\pi;\Pi)$ with $End_\pi(\Pi)$ with values in 
${Ext^2_{\mathbb{Z}[\pi]}}(\Pi,\Pi\otimes\Pi)$.
We may also compose with the natural homomorphisms from 
$\Pi\otimes\Pi$ to $\Pi\odot\Pi$ and $\Gamma_W(\Pi)$ to get pairings
with values in $H^4(L;\Pi\odot\Pi)$ and $H^4(L;\Gamma_W(\Pi))$.

Since $h_\phi^*(\xi\cup\xi')=h_\phi^*\xi\cup{h_\phi^*}\xi'$
we have also
\[
f_\phi(\xi\cup\xi')=
(c_L^*\delta_\phi(p_L^*\xi'))\cup\xi
+(c_L^*\delta_\phi(p_L^*\xi))\cup\xi'
\]
for all $\xi,\xi'\in{H^2}(L;\Pi)$.
In particular, if $\xi\in{H^2(\pi;\Pi)}$ then
$f_\phi(\xi\cup\xi')=0$, and so $f_\phi(c^2_{\pi,w}(\xi))=0$.
Since $c^2_{\pi,w}$ is surjective and the quotient of 
$Ext^2_{\mathbb{Z}[\pi]}(\Pi,\Gamma_W(\Pi))$ 
by the image of $Ext^2_{\mathbb{Z}[\pi]}(\Pi,\Pi\otimes\Pi)$
has exponent 2, by Lemma 12,
it follows that $2f_\phi=0$ on $Ext^2_{\mathbb{Z}[\pi]}(\Pi,\Gamma_W(\Pi))$.

On passing to $\widetilde{L}\simeq{K(\Pi,2)}$ we find that
\[p_L^*(\xi\cup\xi')(\gamma_\Pi(x))=p_L^*\xi(x)\odot{p_L^*}\xi'(x)\] 
for all $\xi,\xi'\in{H^2(L;\Pi)}$ and $x\in\Pi$.
(To see this, note that the inclusion of $x$ determines 
a map from $CP^\infty$ to $K(\Pi,2)$,
since $[CP^\infty,K(\Pi,2)]=Hom(\mathbb{Z},\Pi)$.
Hence we may use naturality of cup products to
reduce to the case when $K(\Pi,2)=CP^\infty$
and $x$ is a generator of $\Pi=\mathbb{Z}$.)
In particular, if $\Xi=\sigma^*id_\Pi\cup\sigma^*id_\Pi$ then
\[p_L^*(\Xi)=2id_{\Gamma_W(\Pi)}\quad\mathrm{ and}\quad 
f_\phi(\Xi)=2(c_L^*\phi)\cup{id_\Pi}=2c^2_{\pi,w}(\phi).\]

If $k=k_2(X)$ for some minimal $PD_4$-complex $X$ with $\pi_1(X)\cong\pi$
then $p_L^*k$ is an isomorphism.
After composition with an automorphism of $\Gamma_W(\Pi)$ 
we may assume that $p_L^*k=id_{\Gamma_W(\Pi)}$,
and so $p_L^*(2k-\Xi)=0$.
Therefore $4(f_\phi(k)-c^2_{\pi,w}(\phi))=2f_\phi(2k-\Xi)=0$.
Since $\mathbb{Z}^w\otimes_\pi\Gamma_W(\Pi)$ is 2-torsion-free
$f_\phi(k)=c^2_{\pi,w}(\phi)$.
Since $c^2_{\pi,w}$ is surjective the orbit of $k$ 
under the action of $E_\pi(L)$
corresponds to an element of 
$Ext^2_{\mathbb{Z}[\pi]}(\Pi,\Pi/2\Pi)\cong{H^2}(\pi;\mathbb{F}_2)$,
and so there are at most $2^\beta$ possibilities.
\end{proof}

If $\pi$ is a $PD_2$-group and $w=w_1(\pi)$ then $L=K\times{CP^\infty}$,
$p_L^*k=z^2\otimes\eta$ and $f_{t\otimes{x}}(k)=2tz\otimes\eta
=c^2_{\pi,w}(t\otimes{x})$.
(However $\mathbb{Z}^w\otimes_\pi\Gamma_W(\Pi)=Z/2Z$ if $w\not=1$.)
The hypotheses also hold if $\pi\cong {Z*_m}$ 
for $m$ even and $w=1$.
(See \S7 and \S9.)

If we could show that $h_\phi$ is the identity on the image of 
$Ext^2_{\mathbb{Z}[\pi]}(\Pi,\Gamma_W(\Pi))$ 
it would follow that $f_\phi$ is additive as a function of $\phi$.
We could then relax the hypothesis on 2-torsion to require 
only that the image of $\Pi\odot_\pi\Pi$ in 
$\mathbb{Z}^w\otimes_\pi\Gamma_W(\Pi)$ be $2$-torsion-free.
(The latter condition holds for all $PD_2$-groups and orientation
characters $w$,
and is easier to check; see Lemma 16 below for the case $\pi=Z*_m$.)

\begin{cor}
If $H^2(\pi;\mathbb{F}_2)=0$, $c^2_{\pi,w}$ is surjective
and $\Pi\odot_\pi\Pi$ is $2$-torsion-free 
there is an unique minimal $PD_4$-complex realizing $(\pi,w)$.
Hence two $PD_4$-complexes $X$ and $Y$ with $\pi_1(X)\cong\pi_1(Y)\cong\pi$ 
are homotopy equivalent if and only if there is an isomorphism 
$\theta:\pi_1(X)\to\pi_1(Y)$ such that $w_1(X)=w_1(Y)\circ\theta$ 
and an isometry of homotopy intersection pairings
$\lambda_X\cong\theta^*\lambda_Y$.
\qed
\end{cor}

We note that we do not yet have explicit invariants that might distinguish
two such minimal $PD_4$-complexes when $\beta>0$.
Does $v_2(X)$ suffice when $\beta=1$?

\section{verifying the torsion condition for $Z*_{m}$}

If $\pi$ is finitely presentable and $c.d.\pi=2$ but $\pi$ 
is not a $PD_2$-group then $H^2(\pi;\mathbb{Z}[\pi])$ 
is not finitely generated \cite{Fa74}.
Whether it must be free abelian remains an open question.
We shall verify this for the groups of most interest to us here.

\begin{lemma}
Let $\pi$ have one end, and be either a semidirect product ${F(s)\rtimes{Z}}$,
a torsion-free one-relator group or the fundamental group of a $3$-manifold 
$M$ with nonempty aspherical boundary.
Then there is a finite $2$-dimensional $K(\pi,1)\!$-complex
and $\Pi=E^2_\pi\mathbb{Z}$ is free abelian.
In particular, $\pi$ is a $2$-dimensional duality group.
\end{lemma}

\begin{proof}
If $\pi=\nu\rtimes{Z}$, where $\nu\cong{F(s)}$ is a nontrivial
finitely generated free group, then $s\geq1$, since $\pi$ has one end.
We may realize $K(\pi,1)$ as a mapping torus of a self-map of $\vee^sS^1$.
This is clearly a finite aspherical 2-complex.
An LHS spectral sequence argument shows that
$\Pi|_\nu=E^2_\pi\mathbb{Z}|_\nu\cong{E^1_\nu\mathbb{Z}}$,
which is free abelian.

If $\pi$ has a one-relator presentation and is torsion-free
the 2-complex associated to the presentation is aspherical 
(and clearly finite).
It is shown in \cite{MT} that one-relator groups are semistable at infinity 
and hence that $\Pi$ is free abelian.

Let $M$ be a $3$-manifold. 
If $\pi=\pi_1(M)$ has one end then
$H_2(\widetilde{M},\partial\widetilde{M};\mathbb{Z})=
H^1(M;\mathbb{Z}[\pi])=0$,
by Poincar\'e duality.
Hence $H_1(\partial\widetilde{M};\mathbb{Z})=0$.
If $\partial{M}$ is a union of aspherical surfaces
it follows that $H_2(\partial\widetilde{M};\mathbb{Z})=0$.
Hence $H_*(\widetilde{M};\mathbb{Z})=0$ for $*>0$ and so $M$ is aspherical.
If moreover $\partial{M}$ is nonempty $M$ retracts onto a finite 2-complex.
The group $\Pi=\overline{H^2(M;\mathbb{Z}[\pi])}$ is free abelian since 
${H^2}(M;\mathbb{Z}[\pi])\cong
{H_1}(\widetilde{M},\partial\widetilde{M};\mathbb{Z})$ 
is the kernel of the augmentation 
$H_0(\partial\widetilde{M};\mathbb{Z})\to{H_0}(\widetilde{M};\mathbb{Z})$.

Since $H^s(\pi;\mathbb{Z}[\pi])=0$ for $s\not=2$ and 
$H^2(\pi;\mathbb{Z}[\pi])$ is torsion-free $\pi$ is a 2-dimensional duality
group \cite{Bi}.
\end{proof}

The class of groups covered by this lemma includes all $PD_2$-groups,
classical knot groups and solvable $HNN$ extensions $Z*_m$ other than $Z$.
Whether every finitely presentable group $\pi$
of cohomological dimension 2 has a finite 2-dimensional $K(\pi,1)$-complex 
and is semistable at infinity remain open questions.

\begin{lemma}
Let $\pi=Z*_m$, $w=1$ and $\Pi=E^2\mathbb{Z}$.
Then $\Pi\odot_\pi\Pi$ is torsion-free.
\end{lemma}

\begin{proof}
The group $\pi=Z*_m$ has a one-relator presentation 
$\langle{a,t}\mid{ta=a^mt}\rangle$
and is also a semidirect product $Z[\frac1m]\rtimes{Z}$.
Let $R=\mathbb{Z}[\pi]$ and $D=\mathbb{Z}[a_n]/(a_{n+1}-a_n^m)$,
where $a_n=t^nat^{-n}$ for $n\in\mathbb{Z}$.
Then $R=\oplus_{n\in\mathbb{Z}}{t^n}D$ is a twisted Laurent extension of 
the commutative domain $D$.

On dualizing the Fox-Lyndon resolution of the augmentation module 
we see that ${H^2}(\pi;\mathbb{Z}[\pi])\cong{R}/(a^m-1, t-\mu_m)R$
and so $\Pi\cong{R}/R(a^m-1,t\mu_m-1)$, where $\mu_m=\Sigma_{i=0}^{i=m-1}a^i$.
Let $E=D/(a^m-1)$ and let $a_{k/m^n}$ be the image of $a_{-n}^k$ in $E$.
Then $E$ is freely generated as an abelian group by $\{a_x\mid{x}\in{J}\}$,
where $J=\{\frac{k}{m^n}\mid0<n,~0\leq{k}<m^{n+1}\}$.
Since $ta_{1-n}^k=a_{-n}^kt$ we have $\Pi\cong\oplus_{n\in\mathbb{Z}}t^nE/\sim$,
where $t^ma_{x}\sim{t^m}a_xt\mu_m=t^{m+1}\mu_ma_{x/m}$.

Therefore $\Pi\odot\Pi\cong\oplus_{m\in\mathbb{Z}}(t^mE\odot{t^mE})/\sim$,
where 
\[t^ma_x\odot{t^ma_y}\sim{t^{m+1}}\mu_ma_{x/m}\odot{t^{m+1}}\mu_ma_{y/m}.\]
Setting $z=y-x$ this gives
\[t^ma_x(1\odot{a_z})\sim{t^{m+1}}a_{x/m}(\mu_m\odot\mu_m{a_{z/m}})
={t^{m+1}}a_{x/m}
(\Sigma_{i,j=0}^{i,j=m-1}a^i(1\odot{a^{j-i}a_{z/m}})).\]
Define a function $f:E\to\Pi\odot\Pi$ by $f(e)=1\odot{e}=e\odot1$ for $e\in{E}$.
Then $f(a_x)=a_xf(a_{m-x})$ for all $x$,
since $a_x\odot1=a_x(1\odot{a_{m-x}})$.
On factoring out the action of $\pi$ we see that 
\[\Pi\odot_\pi\Pi\cong{E}/
(a_z-a_{m-z},a_z-m(\Sigma_{k=0}^{k=m-1}a^ka_{z/m})~\forall{z\in{J}}).\]
(In simplifying the double sum we may set $k=j-i$ for $j\geq{i}$ and 
$k=j+m-i$ otherwise, since $a^ma_{z/m}=a_{z/m}$ for all $z$.)
Thus $\Pi\odot_\pi\Pi$ is a direct limit of free abelian groups and so is
torsion-free.
\end{proof}

If moreover $\mathbb{Z}\otimes_\pi\Pi/2\Pi=H^2(\pi;\mathbb{F}_2)=0$ 
then $\Pi\odot_\pi\Pi\cong\mathbb{Z}\otimes_\pi\Gamma_W(\Pi)$.
Thus if $\pi=Z*_m$ with $m$ even 
$\mathbb{Z}\otimes_\pi\Gamma_W(\Pi)$ is torsion-free.
This group is also torsion-free for $Z*_1=Z^2$; 
does this hold for all $m$?

\section{applications to 2-knots}

Suppose that $\pi$ is either the fundamental group of a finite graph of groups,
with all vertex groups $Z$, or is square root closed accessible,
or is a classical knot group.
(This includes all $PD_2$-groups, semidirect products
$F(s)\rtimes{Z}$ and the solvable groups $Z*_m$.)
Then $L_5(\pi,w)$ acts trivially on
the $s$-cobordism structure set $S^s_{TOP}(M)$ 
and the surgery obstruction map $\sigma_4(M):[M,G/TOP]\to{L_4(\pi,w)}$
is onto, for any closed 4-manifold $M$ realizing $(\pi,w)$.
(See Lemma 6.9 and Theorem 17.8 of \cite{Hi}.)
Thus there are finitely many $s$-cobordism classes within each homotopy type
of such manifolds.

In particular, $Z*_m$ has such a graph-of-groups structure and is solvable, 
so the 5-dimensional TOP $s$-cobordism theorem holds.
Thus if $m$ is even the closed orientable 4-manifold $M$
with $\pi_1(M)\cong{Z*_m}$ and $\chi(M)=0$ is unique up to homeomorphism.
If $m=1$ there are two such homeomorphism types,
distinguished by the second Wu class $v_2(M)$.

Let $\pi$ be a finitely presentable group with $c.d.\pi=2$.
If $H_1(\pi;\mathbb{Z})=\pi/\pi'\cong{Z}$ and $H_2(\pi;\mathbb{Z})=0$ 
then $\mathrm{def}(\pi)=1$, by Theorem 2.8 of \cite{Hi}
If moreover $\pi$ is the normal closure of a single element
then $\pi$ is the group of a 2-knot $K:S^2\to{S^4}$.
(If the Whitehead Conjecture is true every knot group of deficiency 1 
has cohomological dimension at most 2.)
Since $\pi$ is torsion-free it is indecomposable, by a theorem of Klyachko.
Hence $\pi$ has one end. 

Let $M=M(K)$ be the closed 4-manifold obtained by surgery on the 2-knot $K$.
Then $\pi_1(M)\cong\pi=\pi{K}$ and $\chi(M(K))=\chi(\pi)=0$,
and so $M$ is a minimal model for $\pi$.
If $\pi=F(s)\rtimes{Z}$ the homotopy type of $M$ is determined by $\pi$,
as explained in \S4 above.
If $K$ is a ribbon 2-knot it is -amphicheiral
and is determined (up to reflection) by its exterior.  
It follows that a fibred ribbon 2-knot is determined up to $s$-concordance 
and reflection by its fundamental group together with the conjugacy class 
of a meridian.
(This class of 2-knots includes all Artin spins of fibred 1-knots.)

A stronger result holds for the group $\pi=Z*_2$.
This is the group of Fox's Example 10, which is a ribbon 2-knot \cite{Fo62}.
In this case $\pi$ determines the homotopy type of $M(K)$, 
by Theorem 13.
Since metabelian knot groups have an unique conjugacy class
of normal generators (up to inversion) 
Fox's Example 10 is the unique 2-knot 
(up to TOP isotopy and reflection) with this group.
This completes the determination of the 2-knots
with torsion-free elementary amenable knot groups.
(The others are the unknot and the Cappell-Shaneson knots.
See Chapter 17.\S6 of \cite{Hi} for more on 2-knots with  $c.d.\pi=2$.)

\end{document}